\documentclass{amsart}
\usepackage{amssymb}
\usepackage{amsthm}
\usepackage{graphics}
\newtheorem*{reftheorem}{Theorem}
\newtheorem{introtheorem}{Theorem}
\newtheorem{theorem}{Theorem}
\newtheorem{corollary}{Corollary}
\newtheorem{introcorollary}{Corollary}

\newtheorem{lemma}{Lemma}

\newtheorem*{refconj}{Conjecture}
\newtheorem{introconj}{Conjecture}
\newtheorem*{introconj1a}{Conjecture 1A}
\newtheorem*{introconj2a}{Conjecture 2A}
\newtheorem*{theorem1}{Theorem 1}

\newtheorem*{ques}{Question}
\numberwithin{equation}{section}
\title{Volumes of Balls in Large Riemannian Manifolds}
\author{Larry Guth}
\address{Department of Mathematics, Stanford, Stanford CA, 94305 USA}
\email{lguth@math.stanford.edu}

\begin{document} 
\begin{abstract}We prove two lower bounds for the volumes of
balls in a Riemannian manifold.  If $(M^n, g)$ is a complete
Riemannian manifold with filling radius at least $R$, then
it contains a ball of radius $R$ and volume at least
$\delta(n) R^n$.  If $(M^n, hyp)$ is a closed hyperbolic manifold
and if $g$ is another metric on $M$ with volume at most
$\delta(n) Vol(M, hyp)$, then the universal cover of $(M,g)$
contains a unit ball with volume greater than the volume of a
unit ball in hyperbolic n-space.
\end{abstract}

\maketitle

Let $(M, g)$ be a Riemannian manifold of dimension n.  Let
$V(R)$ denote the largest volume of any metric ball of radius $R$
in $(M, g)$.  In \cite{G1}, Gromov made a number of conjectures
relating the function $V(R)$ to other geometric invariants of
$(M,g)$.  The spirit of these conjectures is that if $(M,g)$ is
``large'', then $V(R)$ should also be large.  In this paper, we
prove some of Gromov's conjectures.

Our first result involves the filling radius of $(M,g)$, defined
in \cite{G2}.  Roughly speaking, the filling radius describes how
``thick'' a Riemannian manifold is.  For example, the standard product
metric on the cylinder $S^1 \times \mathbb{R}^{n-1}$ has filling radius
$\pi/2$, and the Euclidean metric on $\mathbb{R}^n$ has infinite
filling radius.

\begin{introtheorem} For each dimension n, there is a number
$\delta(n) > 0$ so that the following estimate holds.  If
$(M^n,g)$ is a complete Riemannian n-manifold with filling radius
at least $R$, then $V(R) \ge \delta(n) R^n$.
\end{introtheorem}

Our second result involves a closed hyperbolic manifold $(M,
hyp)$ equipped with an auxiliary metric $g$.  Slightly
paradoxically, if the manifold $(M, g)$ is small, then its
universal cover $(\tilde M, \tilde g)$ tends to be large.  For
example, if we look at the universal cover of $(M, \lambda^2
hyp)$, we get the space form with constant curvature $-
\lambda^{-2}$.  As $\lambda$ decreases, the strength of the
curvature increases, which increases the volumes of balls.  Our
second theorem gives a large ball in the universal cover $(\tilde
M, \tilde g)$ provided that the volume of $(M, g)$ is
sufficiently small.

\begin{introtheorem} For each dimension n, there is a number
$\delta(n) > 0$ so that the following estimate holds.  Suppose
that $(M^n, hyp)$ is a closed hyperbolic n-manifold and that $g$
is another metric on $M$, and suppose that $Vol(M,g) < \delta(n)
Vol(M,hyp)$.  Let $(\tilde M, \tilde g)$ denote the universal
cover of $M$ with the metric induced from $g$.  Then there is
a point $p \in \tilde M$ so that the unit ball around $p$ in 
$(\tilde M, \tilde g)$ has a larger volume than the unit ball
in hyperbolic n-space.  In other words, the following inequality holds.

$$V_{(\tilde M, \tilde g)}(1) > V_{\mathbb{H}^n}(1).$$
\end{introtheorem}

We spend most of this introduction giving a context for these two
results.  At the end, we give a quick overview of the proof.

Many readers may not be familiar with the filling radius.  Before
looking at its definition, we give some corollaries of Theorem 1
using more common vocabulary.

\begin{introcorollary} Let $(M^n, g)$ be a closed Riemannian
manifold.  Suppose that there is a degree 1 map from $(M^n, g)$
to the unit n-sphere with Lipschitz constant 1.  Then $V(R) \ge
\delta(n) R^n$ for all $R \le 1$.
\end{introcorollary}

\begin{introcorollary} (Systolic inequality) Let $(M^n, g)$ be a
closed aspherical Riemannian manifold.  Suppose that the shortest
non-contractible curve in $(M^n, g)$ has length at least $S$.  Then $V(S)
\ge \delta(n) S^n$.
\end{introcorollary}

\begin{introcorollary} Let $(M^n, g)$ be a closed aspherical
Riemannian manifold, and let $V(R)$ measure the volumes of balls
in the universal cover $(\tilde M, \tilde g)$.  Then $V(R) \ge
\delta(n) R^n$ for all $R$.
\end{introcorollary}

\centerline{\bf Background on filling radius}

We next review the definition of filling radius and some main
facts about it.  For much more information, see \cite{G2}. The
filling radius of a Riemannian manifold is defined by analogy
with an invariant for submanifolds of Euclidean space.  Let $M^n
\subset \mathbb{R}^N$ be a closed submanifold of Euclidean space. 
By a filling of $M$, we mean an (n+1)-chain $C$ with boundary
$M$.  (If $M$ is oriented the standard convention is to use a
chain with integral coefficients, and if $M$ is not oriented the
standard convention is to use a chain with mod 2 coefficients.) 
The filling radius of $M$ is the smallest number $R$ so that $M$
can be filled inside of its $R$-neighborhood.  For example, the
filling radius of an ellipse is its smallest principal axis.  The
main result about filling radius in Euclidean space is the
following estimate.

\begin{reftheorem}(Federer-Fleming, Michael-Simon) If $M^n
\subset \mathbb{R}^N$ is a closed submanifold, then its filling
radius is bounded in terms of its volume by the following
formula.

$$Fill Rad(M) \le C_n Vol(M)^{1/n}.$$
\end{reftheorem}

In \cite{FF}, Federer and Fleming gave a direct construction to
prove that the filling radius of $M$ is bounded by $C_N
Vol(M)^{1/n}$.  Their result is slightly weaker than the result
above because their constant $C_N$ depends on the ambient
dimension $N$.  In \cite{MS}, Michael and Simon proved an
isoperimetric inequality for minimal surfaces that implies the
above theorem.  In \cite{BS}, Bombieri and Simon established the
sharp constant, which occurs when $M$ is a round sphere.

In \cite{G2}, Gromov defined an analogous filling radius for a
closed Riemannian manifold $(M,g)$.  The key observation is that
$(M,g)$ admits a canonical isometric embedding into the Banach
space $L^{\infty}(M)$.  The embedding, which goes back to
Kuratowski, sends the point $x \in M$ to the function $dist_x$
defined by $dist_x(z) = dist(x,z)$.  This embedding depends only
on the metric $g$, and it is isometric in the strong sense that
$|dist_x - dist_y|_{L^\infty} = dist(x,y)$.  Gromov defined the
filling radius of $(M,g)$ to be the infimal $R$ so that the image
of $M$ in $L^\infty$ can be filled inside its $R$-neighborhood. 
There is a similar definition for any complete Riemannian
manifold. 

This definition may seem abstract at first, but \cite{G2}
contains a number of estimates that make it a useful tool.  Here
are a few facts to give a flavor for it.  The filling radius of
the Euclidean metric on $\mathbb{R}^n$ is infinite, but the
filling radius of the standard product metric on $S^n \times
\mathbb{R}^q$ is finite for $n \ge 1$.  If there is a degree 1
map from $(M,g)$ to $(N,h)$ with Lipschitz constant 1, then the
filling radius of $M$ is at least the filling radius of $N$.  The
most important result about filling radius is an analogue of the
Euclidean estimate above.

\begin{reftheorem} (Gromov, \cite{G2}) If $(M^n, g)$ is a
complete Riemannian manifold of dimension $n$, then its filling
radius can be bounded in terms of its volume by the following
formula.

$$Fill Rad(M,g) \le C_n Vol(M,g)^{1/n}.$$

\end{reftheorem}

In the Euclidean setting, Gromov found a stronger version of the
filling radius estimate.  (This result appears near the end
of section F of appendix 1 of \cite{G2}.)

\begin{reftheorem} (Gromov, \cite{G2}, Local volume estimate) If
$M^n \subset \mathbb{R}^N$ has filling radius at least $R$, then
there is some point $x \in \mathbb{R}^N$ so that the volume of $M
\cap B(x,R)$ is at least $c_n R^n$.
\end{reftheorem}

In addition to giving a lower bound for the total volume of $M$,
this result also controls the way the volume is distributed.  It
rules out the possibility that $M$ could have a large total
volume distributed in a diffuse way.  After explaining this
result, Gromov raised the question whether this local volume
estimate has an analogue for Riemannian manifolds.  Our first
theorem answers this question in the affirmative.

\vskip5pt

\centerline{\bf Relation to entropy estimates}

Our second theorem is related to an inequality of Besson,
Courtois, and Gallot.  Their inequality bounds the entropy of a
Riemannian manifold, which is a way of describing the asymptotic
behavior of the volumes of large balls.  In our language, their
result goes as follows.

\begin{reftheorem} (Besson, Courtois, and Gallot, \cite{BCG}) Let
$(M^n, hyp)$ be a closed hyperbolic manifold, and let
$g$ be another metric on $M$ with $Vol(M,g) < Vol
(M,hyp)$.  Then there is some constant $R_0$ (depending on $g$),
so that for every radius $R > R_0$, the following inequality
holds.

$$V_{(\tilde M, \tilde g)} (R) > V_{\mathbb{H}^n}(R).$$

\end{reftheorem}

Our theorem is not as sharp as the theorem of Besson, Courtois,
and Gallot.  The sharp constant in their theorem is a major
achievement.  (The result was previously proven by Gromov with a
non-sharp constant in \cite{G5}.)  To complement their theorem,
it would be nice to estimate the value of $R_0$.  It even looks
plausible that the theorem remains true with $R_0 = 0$!  We
discuss this possibility more below.

Our second theorem can be looked at as a step towards estimating
$R_0$.  According to Theorem 2, the stronger hypothesis $Vol(M,g)
< \delta(n) Vol(M,hyp)$ implies the conclusion $V_{(\tilde M,
\tilde g)}(1) > V_{\mathbb{H}^n}(1)$.  Our method can be modified
to give a similar estimate for balls of radius $R$, but as $R$
moves away from 1, the hypothesis gets stronger (and so the
result gets weaker).  For each $R$, there is a constant,
$\delta(n, R) > 0$ so that $Vol(M,g) < \delta(n, R) Vol(M,hyp)$
implies $V_{(\tilde M, \tilde g)}(R) > V_{\mathbb{H}^n}(R)$.  As
$R$ goes to infinity, the constant $\delta(n, R)$ falls off
exponentially or faster.  As $R$ increases, the methods in this
paper become less effective, whereas the methods in \cite{G5} and
\cite{BCG} are only effective asymptotically for very large R. 
Perhaps there is some way to combine the approaches to get a
uniform estimate for $R \ge 1$.

\vskip5pt

\centerline{\bf Questions about the sharp constants}

It would be interesting to know the sharp constants in Theorems 1
and 2.  In \cite{G1}, Gromov made the following sharp conjecture.

\begin{refconj} (Gromov) Let $(M^n, g)$ be a complete
Riemannian manifold with infinite filling radius.  Let $\omega_n$
be the volume of the unit n-ball.  Then $V(R) \ge \omega_n R^n$.
\end{refconj}

Imitating Gromov, we mention conjectures about the sharp
constants in theorems 1 and 2.  I find the conjectures
intriguing, but the evidence for them is not very strong.

\begin{introconj} Let $(M^n, g)$ be a complete Riemannian
manifold, and let $(S^n, g_0)$ be a round sphere.  We choose the
radius of the round sphere so that the filling radius of $(M^n,
g)$ is equal to that of $(S^n, g_0)$.  Then the following
inequality holds for all $R$.

$$V_{(M^n, g)}(R) \ge V_{(S^n, g_0)}(R).$$

\end{introconj}

\begin{introconj} Let $(M^n, hyp)$ be a closed hyperbolic
manifold, and let $g$ be another metric on $M$ with $Vol(M,g) <
Vol(M, hyp)$.  Then the following inequality holds for all $R$.

$$V_{(\tilde M, \tilde g)} (R) > V_{\mathbb{H}^n}(R).$$
\end{introconj}

As Gromov pointed out in \cite{G1}, estimates of $V(R)$ for small
$R$ are related to scalar curvature.  If $p$ is a point in a
Riemannian manifold $(M^n, g)$, then for small radii $R$, the
volume of the ball $B(p,R)$ is equal to $\omega_n [R^n - (6
n)^{-1} Sc(p) R^{n+2} + o(R^{n+2})]$.  Therefore, the above
conjectures imply several important open conjectures about scalar
curvature.  The scalar curvature conjectures are known in some
special cases, giving modest evidence in favor of the conjectures
above.

\begin{refconj} (Gromov, \cite{G1}) If $(M,g)$ is a complete
Riemannian manifold with infinite filling radius, then $inf_M Sc
\le 0$.  Therefore, if $N$ is a closed aspherical manifold, then
$N$ does not admit a metric of positive scalar
curvature.
\end{refconj}

The last part of the conjecture is known for many particular
aspherical manifolds, but not for all of them.

\begin{introconj1a} (Gromov, \cite{G6}) If $(M^n,g)$ is a
complete Riemannian manifold with scalar curvature at least that
of the unit n-sphere, then the filling radius of $(M,g)$ is at most
that of the unit n-sphere.
\end{introconj1a}

The filling radius of the unit n-sphere was computed by Katz in \cite{K}.  
Gromov and Lawson proved this conjecture with a non-sharp
constant for $n=3$, at least if $H_1(M) = 0$.  (See page 129 of \cite{G2}
and \cite{GL}.)  In higher dimensions, it is unknown.

\begin{introconj2a} (Schoen, \cite{S}) If $(M^n, hyp)$ is a
closed hyperbolic n-manifold and $g$ is another metric on $M$
with scalar curvature at least the scalar curvature of hyperbolic 
n-space, then the volume of $(M,g)$ is
at least the volume of $(M, hyp)$.
\end{introconj2a}

The conjecture is true for $n=2$ by the Gauss-Bonnet theorem. 
For $n=3$, it follows as a corollary of Perelman's proof of
geometrization.  For $n=4$, the conjecture is unknown for
hyperbolic manifolds, but LeBrun \cite{LB} proved a completely
analogous result for certain other 4-manifolds, including the
product of two hyperbolic surfaces.  LeBrun's proof uses
Seiberg-Witten theory.  In higher dimensions, the conjecture is
open.

\vskip5pt
\centerline{\bf Quick summary of the paper}

The main idea of the proof - which is due to Gromov - is to cover
$(M, g)$ with balls, and look at the map from $M$ to the nerve of
the cover.  This technique works well if the multiplicity of the
covering is bounded, because then the dimension of the nerve and
the Lipschitz constant of the map are both under control.  If the
Ricci curvature of $g$ is bounded below, then one can construct a
cover by unit balls with bounded multiplicity, as in \cite{G5}. 
In this paper, we work with only an upper bound on volumes of
balls, and I don't see how to construct a useful cover with
bounded multiplicity.  Instead, we construct a cover that has
bounded multiplicity at most points and push Gromov's ideas to
work on this cover.  For sections 1-5, we assume that $(M,g)$ is
closed.  In section 6, we explain some minor technicalities to
deal with the open case.

Throughout the paper, we use the following notation.  If $U$ is a
region in $(M,g)$, then we denote the volume of $U$ by $|U|$.  If
$z$ is a cycle or chain, then we denote the mass of $z$ by $|z|$. 
If $B$ is shorthand for a ball $B(p, r)$ - the ball around $p$ of
radius $r$ - then $2 B$ is shorthand for $B(p, 2r)$.  Unless
otherwise noted, the constants that appear depend only on the
dimension $n$.

\section{Good balls and good covers}

In this section, we construct a covering of $(M,g)$ by balls with
certain good properties.  All the material in this section is due
to Gromov and appears in sections 5 and 6 of \cite{G2}.

Let $B(p, R) \subset M$ denote the ball around p of radius R. 
We say that the ball $B(p,R)$ is a good
ball if it satisfies the following conditions.

A. Reasonable growth: $|B(p, 100R)| \le 10^{4(n+3)} |B(p, 100^{-1} R)|$.

B. Volume bound: $|B(p,R)| \le 10^{2n+6} V(1) R^{n+3}$.

C. Small radius: $R \le (1/100)$.

The exact constants here are not important.  Notice that in
Euclidean space we would have $|B(p, 100R)| = 10^{4n} |B(p,
100^{-1} R)|$.  The reasonable growth condition relaxes this
bound by replacing $10^{4n}$ with $10^{4(n+3)}$.  In the volume
bound, in Euclidean space we would have $|B(p,R)| = V(1) R^n$. 
For small R, our bound $10^{2n+6} V(1) R^{n+3}$ is much stronger
than the Euclidean bound.  So good balls with small radii have
tiny volumes.  The only crucial point is that $n+3 > n$. 
The other constants were chosen to guarantee the following lemma.

\begin{lemma} Let $(M^n, g)$ be a complete Riemannian n-manifold,
and let $p$ be any point in $M$.  Then there is a radius $R$ so
that $B(p,R)$ is a good ball.
\end{lemma}

\proof We define the density of a ball $B(p,R)$ to be the ratio
$|B(p,R)| / R^n$.  If $B(p,R)$ does not have reasonable growth,
then the density falls off according to the following inequality.

$$\textrm{Density}[B(p, 100^{-1} R)] < 10^{-12}
\textrm{Density}[B(p, 100 R)]. \eqno{(*)}$$

We consider the sequence of balls around $p$ with radii $10^{-2},
10^{-6}, 10^{-10},$ and so on.  We first claim that one of these
balls has reasonable growth.  If the claim is false, then we can
repeatedly use inequality $(*)$ to show that the density of $B(p,
10^{-4s})$ is at most $10^{-12 s} V(1)$.  Since $(M,g)$ is a
Riemannian manifold, the density of $B(p, \epsilon)$ approaches
the volume of the unit n-ball as $\epsilon$ goes to zero.  This
contradiction shows that one of the balls in our list has
reasonable growth.

Now we define $s$ so that $B(p, 10^{-4s-2})$ is the first ball in
the list with reasonable growth.  Applying inequality $(*)$ to the
previous balls $B(p, 10^{-2}), ..., B(p, 10^{-4s+2})$, we
conclude that the density of $B(p, 10^{-4s})$ is at most $10^{-12 s}
V(1)$.  In other words, $|B(p, 10^{-4s})| \le 10^{-12 s} V(1)
[10^{-4s}]^n$.  The ball $B(p, 10^{-4s-2})$ is contained in $B(p,
10^{-4s})$, so its volume obeys the following bound.

$|B(p, 10^{-4s-2})| \le 10^{-4s(n+3)}
V(1) \le 10^{2n+6} [10^{-4s-2}]^{n+3} V(1)$.  

\noindent In other words, the ball $B(p, 10^{-4s-2})$ obeys
condition B.  It has radius $10^{-4s-2} \le 10^{-2}$, and so it
obeys condition C.  Therefore, it is a good ball. \endproof

Because of Lemma 1, we can cover $(M,g)$ with good balls.  We now
use the Vitali covering lemma to choose a convenient sub-covering
with some control of the overlaps.  More precisely, we call an
open cover $\{ B_i \}$ good if it obeys the following properties.

1. Each open set $B_i$ is a good ball.

2. The concentric balls $(1/2) B_i$ cover $M$.

3. The concentric balls $(1/6) B_i$ are disjoint.

(Recall that if $B_i$ is short-hand for $B(p_i, r_i)$, then $(1/2)
B_i$ is short-hand for $B(p_i, (1/2) r_i)$.)

\begin{lemma} If $(M^n, g)$ is a closed Riemannian manifold,
then it has a good cover.
\end{lemma}

\proof This follows immediately from the Vitali covering
lemma.  For each point $p \in M$, pick a good ball $B(p)$.  Then
look at the set of balls $\{ (1/6) B(p) \}_{p \in M}$.  These
balls cover $M$.  Applying the Vitali covering lemma to this set of
balls finishes the proof. \endproof

We now fix a good cover for our manifold $(M,g)$, which we will
use for the rest of the paper.  Our next goal is to control the
amount of overlap between different balls in the cover.  It would
be convenient if we could prove a bound on the multiplicity of
the cover.  I don't see how to prove such a bound, and it may
well be false.  We will prove a weaker estimate in the next
section, bounding the volume of the set where the multiplicity is
high.  We begin with an estimate that controls the number of
balls of roughly equal radius which meet a given ball.

\begin{lemma} If $s<1$, and we look at any ball $B(s)$ of radius
$s$, not necessarily in our cover, then the number of balls $B_i$
from our cover, with radius in the range $(1/2) s \le r_i \le
2s$, intersecting $B(s)$, is less than $C$.
\end{lemma}

\proof Let $\{ B_i \}$ be the set of balls in our cover that
intersect $B(s)$ and have radii in the indicated range.  We
number them so that $B_1$ has the smallest volume.  Now, all of
the balls are contained in $B(5s)$, and $B(5s)$ is contained in
the ball $20 B_1$.  On the other hand, all the $(1/6) B_i$ are
disjoint.  So we have $\sum_i |(1/6) B_i| < |20 B_1|$.  Because
of the locally bounded growth of good balls, we have $\sum |B_i|
< C |B_1|$.  But since $B_1$ has the smallest volume of all the
balls, we see that the number of balls is at most $C$.
\endproof

\section{The volume of the high-multiplicity set}

Let $m(x)$ be the multiplicity function of the cover.  In other
words $m(x)$ is defined to be the number of balls in our cover
that contain the point $x$.  Let $M(\lambda)$ be the set of
points where the multiplicity is at least $\lambda$.  We won't be
able to prove an upper bound for the multiplicity of our
covering, but we partly make up for that by bounding the size of
the set $M(\lambda)$ for large $\lambda$.

We will prove a bound for the total volume $|M(\lambda)|$, but
this bound by itself is not strong enough to prove our theorems. 
We also need bounds on the size of $M(\lambda) \cap B_i$ for
balls $B_i$ in our cover.  For any open set $U \subset M$, we
define $M_U(\lambda) = M(\lambda) \cap U$ to be the set of points
in $U$ with multiplicity at least $\lambda$.  We write $N_w(U)$
to denote the $w$-neighborhood of $U$: the set of points $y \in M$
with $dist(y, U) < w$.

\begin{lemma} There are constants $\alpha(n), \gamma(n)$,
depending only on $n$, that make the following estimate hold. 
For any open set $U \subset M$, and any $w < (1/100)$,

$$|M_U(\gamma \log (1 / w) + \lambda)| \le e^{- \alpha \lambda}
|N_w(U)|.$$

Taking $U = M$, it follows that $|M(\lambda)| < C e^{- \alpha
\lambda} |M|$.

If $B$ is a good ball in our cover with radius $r$, then we have
the following estimate.

$$|M_B(\gamma \log (1/r) + \lambda)| \le e^{-\alpha \lambda}
|B|.$$
\end{lemma}

According to this lemma, it may happen that every point in $B$
has multiplicity $\gamma(n) \log (1 / r)$.  However, the set of
points in $B$ with multiplicity much higher than $\gamma(n) \log
(1/r)$ constitutes only a small fraction of $B$.

\proof Let $\{ B_i \}_{i \in I}$ denote the subset of balls in
our cover that intersect $U$.  We divide this set of balls into
layers, using the Vitali covering construction to choose each
layer.

To choose $Layer(1)$: Take the largest ball in the set.  (More
precisely, take the ball with the largest radius.)  Then take the
next largest ball disjoint from it.  Then take the largest ball
disjoint from the two already chosen... When there are no more
balls left, stop.

To choose $Layer(2)$: Examine all the balls that are not part of
$Layer(1)$.  Take the largest ball available.  Then take the
largest remaining ball disjoint from this one...

To choose $Layer(d)$: Examine all the balls that are not part of
any previous layer.  Take the largest ball available...

In this way, our set of balls $\{ B_i \}_{i \in I}$ is divided
into layers.  Each ball in the set belongs to exactly one layer. 
Each layer consists of disjoint balls.  We define $L(d)$ to be
the union of all the balls in $Layer(d)$.  We call $Layer(1)$ the
top layer, and if $d_2 > d_1$, we say that $Layer(d_2)$ is lower
than $Layer(d_1)$.  Now for each layer, we define a subset
$Core(d) \subset L(d)$ which intersects only a bounded number of
balls from lower layers.

To define the core, we first introduce a partial ordering on the
balls in a given layer $Layer(d)$.  If $B_i, B_j \in Layer(d)$,
we say that $B_i < B_j$ if there is some ball $B_k$ in a lower
layer which meets both $B_i$ and $B_j$, and if the radii obey the
inequalities $2 r_i \le r_k \le r_j$.  We consider the minimal
partial order that is generated by these relations.  In other
words, we say that $B_i < B_j$ if and only if there is a chain of
balls $B_{l_1}, ..., B_{l_m}$ in $Layer(d)$ and a chain of balls
$B_{k_1}, ..., B_{k_{m+1}}$ in lower layers so that $B_i$ meets
$B_{k_1}$ which meets $B_{l_1}$ which meets $B_{k_2}$ ... which
meets $B_{l_m}$ which meets $B_{k_{m+1}}$ which meets $B_j$ and
so that the radii obey $2^{m+1} r_i \le 2^m r_{k_1} \le 2^m
r_{l_1} \le 2^{m-1} r_{k_2} \le ... \le 2 r_{l_m}
\le r_{k_{m+1}} \le r_j$.  The following diagram illustrates the
overlapping balls in case $m=1$.  The balls drawn in solid lines
belong to $Layer(d)$ and those in dotted lines belong to lower
layers.  The smallest ball on the left is $B_i$, and the largest
ball on the right is $B_j$.

\includegraphics{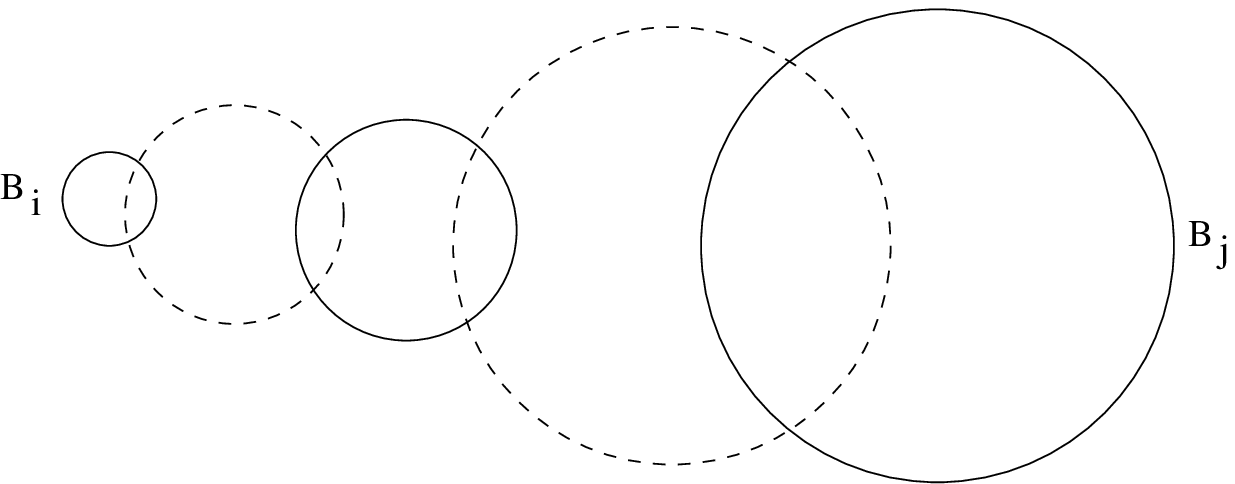}

We define $Max(d) \subset Layer(d)$ to be the maximal elements of
this partial ordering.  (A ball $B_i$ is maximal if there is no
other ball $B_j$ with $B_i < B_j$.)  For any ball $B_i$, we
define the core of $B_i$ to be the concentric ball $(1/10) B_i$. 
We define the core of $Layer(d)$ to be the union of the cores of
all the maximal balls in $Layer(d)$.

$$Core(d) = \cup_{B_i \in Max(d)} \frac{1}{10} B_i.$$

(The balls in $Layer(d)$ are disjoint, so the core is a union of
disjoint balls.)

By looking at maximal balls, we buy the following inequality. 
Let $B_i$ be a maximal ball in $Layer(d)$ and let $x \in (1/10)
B_i$ be a point in $Core(d)$.  Suppose that $x$ also lies in a
ball $B_k$ from a lower layer.  Then the radius $r_k$ is pinched
in the range $(1/15) r_i \le r_k \le 2 r_i$.  The upper bound
depends on the maximality of $B_i$.  Suppose that $r_k > 2 r_i$. 
Since the ball $B_k$ was not selected to join the layer
$Layer(d)$, there must be some larger ball $B_j$ in $Layer(d)$
intersecting $B_k$.  But then it follows that $B_i < B_j$,
contradicting the assumption that $B_i$ is maximal.  The lower
bound for $r_k$ does not depend on maximality.  We know that the
concentric balls $(1/6) B_i$ and $(1/6) B_k$ are disjoint.  In
particular, that implies that the center of $B_k$ lies outside of
$(1/6) B_i$.  Now if $r_k < (1/15) r_i$, then the ball $B_k$ lies
outside of $(1/10) B_i$.  On the other hand, $x$ lies in $(1/10)
B_i$, and we get a contradiction.

According to Lemma 3, the number of balls $B_k$ containing $x$
with radii in the range $(1/15) r_i \le r_k \le 2 r_i$ is bounded
by a dimensional constant $\eta(n)$.  Therefore, the number of
balls $B_i$ so that $x \in B_i$ and so that $B_i \in
Layer(\lambda)$ with $\lambda \ge d$ is at most $\eta(n)$.  Less
formally, this estimate says that the core of $Layer(d)$ is
well-insulated from the balls in lower layers.

The next main point is that $Core(d)$ contains a substantial
fraction of the volume of $L(d)$.  We first claim that $L(d)
\subset \cup_{B_j \in Max(d)} 10 B_j$.  Suppose that $B_i$ is any
ball in $Layer(d)$.  If $B_i$
is itself maximal, then it is contained in the union $\cup_{B_j
\in Max(d)} 10 B_j$.  If $B_i$ is not maximal, then there is some
chain of overlapping balls $B_i, B_{k_1}, B_{l_1}, ..., B_{k_m},
B_{l_m}, B_{k_{m+1}}, B_j$, where $B_j$ is maximal, and the radii
obey $2^{m+1} r_i \le 2^m r_{k_1} \le 2^m r_{l_1} \le ... \le 2
r_{k_m} \le 2 r_{l_m} \le r_{k_{m+1}} \le r_j$.  Because the
balls overlap, $B_i$ is contained in the ball with the same center as
$B_j$ and with radius $R = r_j + 2 r_{k_{m+1}} + 2 r_{l_m} + 2
r_{k_m} + ... + 2 r_{l_1} + 2 r_{k_1} + 2 r_i$.  According to our
bounds for the radii of the balls, $R \le r_j + 2 r_j + 4 \cdot 2^{-1}
r_j + 4 \cdot 2^{-2} r_j + 4 \cdot 2^{-3} r_j + ... \le 7 r_j$.  Therefore,
$B_i \subset 10 B_j$.  Since $B_j$ is maximal, $B_i \subset
\cup_{B_j \in Max(d)} 10 B_j$.  Therefore, $L(d) \subset
\cup_{B_j \in Max(d)} 10 B_j$.  Using this inclusion and the
reasonable growth estimates for good balls, we can estimate the
volume of $Core(d)$.

$$|L(d)| \le |\cup_{B_j \in Max(d)} 10 B_j| \le \sum_{B_j \in
Max(d)} | 10 B_j| \le C \sum_{B_j \in Max(d)} |\frac{1}{10} B_j|
\le C |Core(d)|.$$

We will use these estimates about the core to prove the
exponential decay of the high-multiplicity set.  We now introduce
some vocabulary that describes how many times a point is
contained in balls from different layers.

$$L^\mu(\lambda) := \{ x | x \in L(d) \textrm{ for at least $\mu$
different values of $d$ in the range $d \ge \lambda$} \}.$$

The sets $L^\mu(\lambda)$ are nested: $L^1(\lambda) \supset
L^2(\lambda) \supset ... $ Because of the construction of the
layers, $\cup_{B_i \in Layer(\lambda)} 3 B_i$ contains $\cup_{d
\ge \lambda} \cup_{B_i \in Layer(d)} B_i = L^1(\lambda)$. 
Because the balls in each layer are disjoint and because of the
reasonable growth bound, we get the following upper bound for
$|L^1(\lambda)|$.

$$|L^1(\lambda)| = |\cup_{d \ge \lambda} \cup_{B_i \in Layer(d)}
B_i| \le |\cup_{B_i \in Layer(\lambda)} 3 B_i| \le $$

$$ \le \sum_{B_i \in
Layer(\lambda)} |3 B_i| \le C \sum_{B_i \in Layer(\lambda)} |B_i|
= C |L(\lambda)|.$$

Now we define a function $F(\lambda)$ which is an average of the
volumes of $L^\mu(\lambda)$.

$$F(\lambda) := \frac{1}{\eta(n)} \sum_{\mu=1}^{\eta(n)}
|L^{\mu}(\lambda)|.$$

Our estimates about the core imply that $F(\lambda)$ decays
exponentially.  We proved that each point $x$ in
$Core(\lambda)$ lies in at most $\eta(n)$ balls from layers
$Layer(d)$ with $d \ge \lambda$.  We know that $Core(\lambda)
\subset L(\lambda)$.  Therefore, we get the following estimate.

$$\sum_{\mu=1}^{\eta(n)} |L^\mu(\lambda)| - |L^\mu(\lambda+1)|
\ge |Core(\lambda)|.$$

Using the formula for $F(\lambda)$, we see that
$F(\lambda) - F(\lambda + 1) \ge (1 / \eta) |Core(\lambda)|$. 
Now we plug in our volume estimates for $|Core(\lambda)|$ and
$|L^1(\lambda)|$.

$$F(\lambda) - F(\lambda + 1) \ge (1 / \eta) |Core(\lambda)| \ge
c |L(\lambda)| \ge c' |L^1(\lambda)| \ge c' F(\lambda).$$

Rearranging the equation, we can deduce the exponential decay of
$F(\lambda)$.

$$F(\lambda + 1) \le (1 - c') F(\lambda). \eqno{(*)}$$

(In this equation $c' > 0$ is a small constant depending only on
$n$.)

To control $F(\lambda)$ for large values of $\lambda$, we combine
equation $(*)$ with some simple estimates for $F(\lambda)$ for
small values of $\lambda$.  To control $F(\lambda)$ for small
$\lambda$, we need one more inequality, which says that the large
balls are put into the top layers.  More precisely, if a ball
$B_i$ of radius $r_i$ belongs to $Layer(d)$, then $d \le \gamma
\log (1/r_i)$ (for a constant $\gamma$ depending only on $n$). 
Because $B_i$ does not belong to any of the first $d-1$ layers,
there must be a larger ball in each of those layers intersecting
$B_i$.  All of the balls in our cover have radius at most
$(1/100)$, so the number of different scales of the radii of
these balls is $\log (1/r_i)$.  For each scale, Lemma 3 tells us
that there are at most $C$ balls of that scale intersecting
$B_i$.  Therefore, we can choose a constant $\gamma(n)$ so that
the total number of larger balls intersecting $B_i$ is at most
$\gamma \log (1/r_i)$, and $d \le \gamma \log(1/r_i)$ as claimed. 
So if $B_i$ belongs to $Layer(d)$ with $d \ge \gamma
\log (1/w)$, then $r_i \le w$.  In other words, all the balls in
layers lower than $\gamma \log(1/w)$ have radius at most $w$. 
Since all of our balls intersect $U$, the union of all the balls
in the lower layers is contained in $N_{2w} (U)$.  This argument gives
us the following estimate for $F(\gamma \log(1/w))$.

$$F(\gamma \log(1/w)) \le |L^1(\gamma \log(1/w))| \le |N_{2w}(U)|.$$

Combining this inequality with the exponential decay $(*)$ we get
the following bound.

$$F(\gamma \log(1/w) + \lambda) \le e^{- \alpha \lambda} |N_{2w}(U)|.$$

Finally, we bound the size of $M_U(\lambda)$ in terms of $F$. 
Since each layer consists of disjoint balls, $M_U(\lambda +
\eta(n)) \subset L^{\eta(n)}(\lambda)$, and so $|M_U(\lambda +
\eta(n))| \le |L^{\eta(n)}(\lambda)| \le F(\lambda)$.  Combining
this observation with our last inequality gives the following.

$$|M_U(\gamma \log(1/w) + \lambda + \eta(n))| \le
e^{-\alpha \lambda} |N_{2w}(U)|.$$

This inequality is equivalent to the one we wanted to prove.

As a special case, we can take $U = M$.  Applying our inequality
with $w = (1/100)$ yields $|M(\lambda)| \le C e^{- \alpha
\lambda} |M|$.

As another special case, we can take $U = B$ for a ball $B$ in
our cover of radius $r$.  Applying our inequality with $w = r$
yields

$$|M_B(\gamma \log(1/r) + \lambda + \eta(n))| \le e^{-\alpha \lambda} |2
B|.$$

Since $B$ is a good ball, $|2 B| < C |B|$, and so we get the
following inequality.

$$|M_B(\gamma \log(1/r) + \lambda + \eta(n))| \le C e^{-\alpha \lambda}
|B|.$$

This inequality is equivalent to the one we wanted to prove.
\endproof

\section{The rectangular nerve}

Gromov had the idea to prove estimates about a Riemannian
manifold $(M,g)$ by covering it with balls and considering the
induced map from $M$ to the nerve of the covering.  I believe
that this idea first appeared in \cite{G5}.  It is also discussed
in section 5.32 of \cite{G4}.  In our case, there is an
added wrinkle because the balls in our covering have a wide range
of radii, and we need to choose a metric on the nerve that reflects
the radii of the balls in the covering.  In order to
accomplish that, we slightly modify the idea of the nerve,
introducing a ``rectangular nerve''.

For each ball $B_i$ of radius $r_i$ in our good cover, define
$\phi_i: M \rightarrow [0, r_i]$ as follows.  Outside of $B_i$,
$\phi_i = 0$.  Let $d$ be the distance from $x \in B_i$ to the
center of $B_i$.  If $d \le (1/2) r_i$, then $\phi(x) = r_i$.  If
$(1/2) r_i \le d \le r_i$, then $\phi(x) = 2( r_i - d)$.  The
Lipschitz constant of $\phi_i$ is 2.

All the $\phi_i$ together map $M$ into the high-dimensional
rectangle $R$ with dimensions $r_1 \times ... \times r_D$. 
Because $(1/2) B_i$ covers $M$, we know that for each $x \in M$,
there is some $\phi_i$ so that $\phi_i(x) = r_i$.  Therefore, the
image of M lies in the union of certain hyperfaces of the
rectangle R: namely those closed hyperfaces that don't contain
$0$.  As in the simplicial world, we define a nerve $N$ which
will be a subcomplex of the rectangle $R$.  An open face $F$ of
$R$ is determined by dividing the dimensions $1 ... D$ into three
categories: $I_0$, $I_1$, and $I_{(0,1)}$.  Then $F$ is given by
the equalities and inequalities $\phi_i = 0$ for $i \in I_0$,
$\phi_i = r_i$ for $i \in I_1$, and $0 < \phi_i < r_i$ for $i \in
I_{(0,1)}$.  We denote $I_+ = I_1 \cup I_{(0,1)}$.  A face $F$ is
contained in the nerve if and only if $\cap_{i \in I_+(F)} B_i \not=
\emptyset$, and if $I_1(F) \not= \emptyset$.  We see that $\phi$
maps $M$ into the rectangular nerve $N$.

Using our bounds for the high-multiplicity set in $M$, we can
bound the volume of the image $\phi(M)$, and also the volume of
$\phi(M)$ contained in certain regions of $N$.  If $F$ is an open
face in $N$, then we define $Star(F)$ to be the union of all open
faces of $N$ which contain $F$ in their closures.  In other
words, $Star(F)$ is the union of $F$ itself together with each
higher-dimensional open face that contains $F$ in its boundary. 
No lower dimensional faces are contained in $Star(F)$.  We let
$d(F)$ denote the dimension of $F$.  Each face $F$ is itself a
rectangle with dimensions $r_1(F) \le ... \le r_{d(F)}(F)$. 

\begin{lemma} There are constants $C(n), \beta(n) > 0$, depending
only on $n$, so that the volume of $\phi(M) \cap Star(F)$ obeys
the following inequality.

$$|\phi(M) \cap Star(F)| < C V(1) r_1(F)^{n+1} e^{- \beta d(F)}.$$

Also $|\phi(M)| < C |M|$.

\end{lemma}

\proof The set $\phi^{-1} [Star(F)]$ is contained in $B_1$, the
ball with radius $r_1 = r_1(F)$.  By the good ball estimate, this
set has volume at most $C V(1) r_1^{n+3}$.  The Lipschitz
constant of $\phi$ at a point $x \in B_1$ is bounded by $2
m(x)^{1/2}$.  According to Lemma 4, the set of points in $B_1$
with multiplicity more than $\gamma \log(r_1^{-1}) + \lambda$ is
bounded by $e^{- \alpha \lambda} |B_1|$.  Adding up the
contributions from the regions of different multiplicity, we see
that $|\phi(B_1)| < C V(1) r_1^{n+2}$.  Now, by Lemma 3, it
follows that $r_1(F) < C e^{- \beta d(F)}$.  Plugging in, we get the
lemma.

To prove the last claim, we apply the same argument, using the
bound for $|M(\lambda)|$ in Lemma 4.
\endproof

We should make some remarks about this inequality.  In our paper,
it turns out to be natural to compare $|\phi(M) \cap Star(F)|$
with $r_1(F)^n$.  According to our inequality, the ratio
$|\phi(M) \cap Star(F)|/r_1(F)^n$ is at most $C V(1) r_1(F) e^{-
\beta d(F)}$.  In other words, the ratio becomes favorable if
$V(1)$ is small, or if $r_1(F)$ is small, or if $d(F)$ is large.

To finish this section, we explain the connections between the
rectangular nerve and the filling radius and simplicial volume of
$M$.

\begin{lemma} If $\phi_*([M]) = 0$ in $N$, then the filling
radius of $(M^n, g)$ is at most 1.
\end{lemma}

\proof Let $A$ be a chain in $N$ filling $\phi(M)$.  The filling
$A$ consists of the following data: an abstract chain $A$ with
$\partial A = M$, together with a map $f: A \rightarrow N$ with
$f |_{\partial A} = \phi$.  By abuse of notation, we identify $M$
with its image in $L^\infty(M)$ under the Kuratowski embedding. 
Using the above data, we will construct a filling of $M$ in its
1-neighborhood.

First, we construct a map $\psi$ from $A$ to $L^\infty(M)$.  We
pick a fine triangulation of $A$, subordinate to the faces of
$N$.  We begin by defining $\psi$ on the vertices of the
triangulation.  Each vertex $v$ lies in some face $F$ of the
nerve $N$.  We pick a ball $B_i$ so that the index $i$ lies in
$I_+(F)$.  Then we define $\psi(v)$ to be the point $p_i$ which
is the center of the ball $B_i$.  Next we define $\psi$ on each
simplex by extending it linearly.

Because our triangulation of $A$ is subordinate to the faces of
$N$, each edge of $A$ joins two points that lie in a common
closed face of $N$.  Therefore, $\psi$ maps the endpoints of any
edge to the centers of two overlapping good balls.  The distance
between the centers is at most $(2/100)$.  Each simplex of $A$ is
mapped to a simplex in $L^\infty(M)$ whose edges have length at
most $2/100$.  Also, each vertex is mapped into $M \subset
L^\infty(M)$. Therefore, the image $\psi(A)$ lies in the $2/100$
neighborhood of $M$.

We are not finished, because the map $\psi$ restricted to $M$ is
not the Kuratowksi embedding.  To finish the proof, we will
homotope $M$ to $\psi(M)$ inside of the $2/100$-neighborhood of
$M$.  Let $p$ be a point of $M$.  Let $\Delta$ be the smallest
simplex of our triangulation that contains $p$, and let $y_1,
..., y_m$ be the vertices of $\Delta$.  The map $\psi$ sends each
$y_i$ to the center $p_i$ of some ball $B_i$ containing $y_i$. 
The map $\psi$ sends $p$ to a point on the linear simplex in
$L^\infty$ spanned by the points $p_i$.  The distance from $p$ to
$\psi(p)$ is at most the largest distance from $p$ to any of the
$p_i$.  Since each triangle is small, we can assume that the
distance from $p$ to $y_i$ is less than $1/100$.  Also, since
each good ball has radius at most $1/100$, the distance from
$y_i$ to $p_i$ is at most $1/100$.  Combining these inequalities,
the distance from $p$ to $\psi(p)$ is at most 2/100.  Therefore,
the Kuratowski embedding can be homotoped to the map $\psi$
inside the $2/100$ neighborhood of $M$ in $L^\infty(M)$.

Combining this homotopy with the chain $\psi(A)$, it follows that
$M$ bounds inside its $2/100$ neighborhood.  In other words, the
filling radius of $(M,g)$ is at most $2/100 < 1$. \endproof

\begin{lemma} If $(M,g)$ is a closed aspherical manifold with
systole at least 1, then there is a map $\psi: N \rightarrow M$
so that the composition $\psi \circ \phi: M \rightarrow M$ is
homotopic to the identity.
\end{lemma}

\proof This proof is essentially the same as Gromov's from
\cite{G5} pages 293-294, which uses the standard simplicial nerve
instead of the rectangular nerve.

We slightly homotope $\phi$ to a map $\phi'$ which is simplicial
with respect to some fine triangulations of $M$ and $N$.  We can
assume that the triangulation of $N$ is subordinate to the faces
of $N$.

We begin by defining the map $\psi$ from $N$ to $M$.  We define
the map one skeleton at a time.  For each vertex $v$ of $N$, we
consider the smallest face $F \supset v$, and we pick an index in
$I_+(F)$.  Then we map $v$ to the center of $B_i$.  Now each edge
$E$ of $N$ joins two vertices lying in the same closed face.  If
the boundary of $E$ is $v_1 \cup v_2$, then it follows that we
have mapped $v_1$ and $v_2$ to two overlapping balls from our
covering.  Since each ball has radius at most $(1/100)$, the
distance between the centers is at most $(2/100)$, and we may
map $E$ to an arc of length at most $(2/100)$.  Now the boundary
of each 2-simplex has been mapped to an arc of length at most
$(6/100)$.  Since the 1-systole of $(M,g)$ is at least 1, the
image curve is contractible, and so we can extend our map to each
2-simplex.  Since $M$ is aspherical, we can then extend the map
to each higher-dimensional simplex.  This completes the
construction of $\psi$.

Next we have to show that $\psi \circ \phi'$ is homotopic to the
identity.  We have to define a map $H$ on $M \times [0,1]$ with
$H(m,0) = \psi \circ \phi'(m)$ and $H(m,1) = m$.  We define $H$
one skeleton at a time.  For each vertex $v$, $\phi'(v)$ is a
vertex of the triangulation of $N$ lying very near to $\phi(v)$. 
We let $F(v)$ denote the smallest face of $N$ containing
$\phi'(v)$.  It may not be the case that $\phi(v)$ lies in
$F(v)$, but at least $\phi(v)$ lies in a face bordering $F(v)$. 
Therefore, $\psi \circ \phi'(v)$ is the center of some ball $B_i$
overlapping some other ball $B_j$ containing $v$.  Since each
ball in our cover has radius at most $(1/100)$, the distance from
$\psi \circ \phi'(v)$ to $v$ is at most $(3/100)$.  We define $H$
on $v \times (0,1)$ by mapping the interval to a curve from $\psi
\circ \phi'(v)$ to $v$, with length at most $(3/100)$.

Next we look at an edge $E$ of the triangulation of $M$.  The map
$\phi'$ either collapses $E$ to a point or maps it onto an edge
of the triangulation of $N$.  Therefore, $\psi \circ \phi'(E)$ is
either a point or an arc of length at most $(2/100)$.  We have
already defined $H$ on the boundary of $E \times (0,1)$.  The
restriction of $H$ to the boundary is a curve of length at most
$(8/100)$ plus the length of $E$.  We can assume the length of
$E$ is at most $(1/100)$.  Since the 1-systole of $(M,g)$ is at
least 1, this curve is contractible, and so we can extend $H$ to
$E \times (0,1)$ for every edge E.

Finally, since $M$ is aspherical, we can extend $H$ to $\Delta
\times (0,1)$ for each 2-simplex $\Delta$ of $M$, and then for
each higher-dimensional simplex.  Therefore, $\psi \circ \phi'$
is homotopic to the identity.  Since $\phi'$ is homotopic to
$\phi$, $\psi \circ \phi$ is homotopic to the identity.  \endproof

\section{Filling cycles in the rectangular nerve}

Using the bounds proved in the last section, we will
now show that if $V(1)$ is sufficiently small, then $\phi_*([M])
= 0$ in the rectangular nerve $N$.

\begin{lemma} For any $\beta > 0$ and any integer $n > 0$, there
is a small positive $\epsilon(\beta, n)$ that makes the following
statement true.  Let $z$ be an n-cycle in the rectangular
complex $N$.  Suppose that $z$ obeys the following estimate.

$$|z \cap Star(F)| < \epsilon r_1(F)^n e^{- \beta d(F)}.$$

Then $[z] = 0$ in $N$.
\end{lemma}

\proof Let $D$ be the dimension of $N$.

We will construct a sequence of homologous cycles $z = z_D \sim
z_{D-1} \sim ... \sim z_n$.  The cycle $z_k$ will be
contained in the k-skeleton of $N$.  Moreover, every cycle will
obey the following estimate, slightly weaker than the estimate
that $z$ obeys.

$$|z_k \cap Star(F)| < 2 \epsilon r_1(F)^n e^{- \beta d(F)}.$$

In particular, for each n-face $F$, the cycle $z_n$ obeys the
following estimate.

$$|z_n \cap F| < 2 \epsilon r_1(F)^n.$$

Our constant $\epsilon$ will be less than $1/2$, so we conclude
that $|z_n \cap F| < |F|$.  Therefore, $z_n$ is homologous to a
cycle lying in the (n-1)-skeleton of $N$, and hence $[z_n] = 0$.

Now we describe the inductive step, getting from $z_{k}$ to
$z_{k-1}$.  Let $F$ be a k-dimensional face of $N$.  Consider
$z_k \cap F$, which defines a relative cycle in $F$, and we
replace it with the minimal relative cycle with the same
boundary.  Performing this surgery on each k-face $F$, we get a new
n-cycle $z_k'$, homologous to $z_k$, and still contained in the
k-skeleton of $N$.

We examine the intersection $z_k' \cap F$, for a k-dimensional
face $F$.  Since $z_k'$ was chosen to minimize volume, it follows
that $|z_k' \cap F| \le |z_k \cap F|$.  By the inductive
hypothesis, $|z_k \cap F| \le 2 \epsilon r_1(F)^n e^{- \beta k}$. 
Using this volume estimate, we can show that $z_k'$ lies near to
the boundary of $F$.  Suppose that $x \in z_k'$ and that the
distance from $x$ to $\partial F$ is $s$.  By the monotonicity
formula, it follows that $\omega_n s^n \le |z_k' \cap F| \le 2
\epsilon r_1(F)^n e^{- \beta k}$.  Rearranging this formula, we
get the following inequality, bounding the distance from any
point in $z_k'$ to the boundary $\partial F$.

$$s / r_1(F) \le [2 \omega_n^{-1} \epsilon
e^{-\beta k}]^{1/n}. \eqno{(1)}$$

Following Gromov in \cite{G2} we define a map that pulls a small
neighborhood of the (k-1)-skeleton of $N$ into the
(k-1)-skeleton.  Our map will be called $R_\delta$, and it
depends on a number $\delta$ in the range $0 < \delta < 1/2$. 
The basic map is a map from an interval $[0, r]$ to itself, which
takes the set $[0, \delta r]$ to $0$, and the set $[r- \delta r,
r]$ to $r$, and linearly stretches the set $[\delta r, r - \delta
r]$ to cover $[0, r]$.  The Lipschitz constant of this map is $(1
- 2 \delta)^{-1}$.  Now we apply this map separately to each
coordinate $\phi_i$ of the big rectangle $R$.  The resulting map
is $R_\delta$.

The map $R_\delta$ has the following nice properties.  It maps
the nerve $N$ into itself.  The map $R_0$ is the identity, and so
each $R_\delta$ is homotopic to the identity (in the space of
self-maps of $N$).  Therefore, the map $R_\delta$ moves any cycle
to a homologous cycle. The preimage $R_\delta^{-1} [Star(F)] =
Star(F)$ for any face $F$.  Since the Lipschitz constant of
$R_\delta$ is $[1 - 2 \delta]^{-1}$, the following estimate holds
for any n-cycle $y$.

$$|R_\delta(y) \cap Star(F)| \le (1-2 \delta)^{-n} |y \cap
Star(F)|.  \eqno{(2)}$$

\noindent Also, for sufficiently big $\delta$, the map $R_\delta$
takes $z_k'$ into the (k-1)-skeleton of $N$.  In particular, if
$\delta \ge [2 \omega_n^{-1} \epsilon e^{-\beta k}]^{1/n}$, then
inequality $(1)$ guarantees that $R_\delta(z_k')$ lies in the
(k-1)-skeleton of $N$.  We define $\delta(k) = [2 \omega_n^{-1}
\epsilon e^{- \beta k}]^{1/n}$, and then we define $z_{k-1} =
R_{\delta(k)} (z_k')$.

We have to check that $z_{k-1}$ obeys the volume estimate in the
inductive hypothesis.  Let $F$ be a face of $N$ with any
dimension.  First we claim that it obeys the following estimate.

$$|z_{k-1} \cap Star(F)| \le \prod_{l=k}^D (1- 2 \delta(l))^{-n}
\epsilon r_1(F)^n e^{- \beta d(F)}.$$

This estimate follows from three observations.  First, by
hypothesis, $|z_D \cap Star(F)| \le \epsilon r_1(F)^n e^{-\beta
d(F)}$.  Second, $|z_k' \cap Star(F)| \le |z_k \cap Star(F)|$. 
Third, by equation number 2, $|z_{k-1} \cap Star(F)| \le (1 - 2
\delta(k))^{-n} |z_k' \cap Star(F)|$.  So to make the induction
work, we have to choose $\epsilon$ sufficiently small that the
following estimate holds.

$$ \prod_{l=n+1}^\infty (1 - 2 \delta(l))^{-n} = \prod_{l = n+
1}^ \infty (1 - 2 [2 \omega_n^{-1} \epsilon e^{-
\beta l}]^{1/n})^{-n} < 2.$$

The product converges because of the exponential decay in the
term $e^{- \beta l}$, and by taking $\epsilon > 0$ sufficiently
small, we can guarantee that it is less than 2.  The value of
$\epsilon$ here depends on $n$ and $\beta$. \endproof

We now have enough ammunition to prove Theorem 1 for closed manifolds.

\begin{theorem} (Closed case) For each dimension n, there is
a number
$\delta(n) > 0$ so that the following estimate holds.  If
$(M^n,g)$ is a closed Riemannian n-manifold with filling radius
greater than $R$, then $V(R) \ge \delta(n) R^n$.
\end{theorem}

\proof By scaling, it suffices to prove the theorem when $R = 1$. 
We consider the map $\phi$ from $M$ to the rectangular
nerve $N$ of a good cover.  According to Lemma 5, the image obeys
the following estimate for each face $F$ of $N$.

$$|\phi(M) \cap Star(F)| < C V(1) r_1(F)^{n+1} e^{-\beta d(F)}.$$

\noindent The constants $C, \beta$ in this equation depend only
on $n$.  Let $\epsilon(\beta, n)$ be the number defined in Lemma
8.  Hence there is some number $\delta(n)$ so that if $V(1) <
\delta$, then we get the following estimate for each face $F$ of
$N$.

$$|\phi(M) \cap Star(F)| < \epsilon r_1(F)^n e^{- \beta d(F)}.$$

According to Lemma 8, this estimate implies that the cycle
$\phi(M)$ is homologous to zero in $N$.  Now, according to Lemma
6, the filling radius of $(M, g)$ is at most 1.
\endproof

We now prove two of the corollaries from the introduction.

\begin{corollary} Let $(M^n, g)$ be a closed Riemannian
manifold.  Suppose that there is a degree non-zero map $F$ from
$(M^n, g)$ to the unit n-sphere with Lipschitz constant 1.  Then
$V(R) \ge
\delta(n) R^n$ for $R \le 1$.
\end{corollary}

\proof In \cite{G2} (page 8), Gromov proved that the filling radius of
$(M,g)$ is at least the filling radius of the unit n-sphere.
Theorem 1 then implies that $V(R) > \delta(n) R^n$ for $R \le 1$.
\endproof

\begin{corollary} (Systolic inequality) Let $(M^n, g)$ be a
closed aspherical Riemannian manifold.  Suppose that the shortest
non-contractible curve in $(M^n, g)$ has length at least $S$.  Then $V(S)
\ge \delta(n) S^n$.
\end{corollary}

\proof In \cite{G2} (section 1), Gromov proved that the filling radius of
$(M,g)$ is at least $S/6$.  According to Theorem 1, $V(S/6) \ge
\delta(n) (S/6)^n$.  For a smaller constant $\delta(n)$, $V(S)
\ge V(S/6) \ge \delta(n) S^n$. \endproof

\section{Estimates for simplicial norms}

In this section, we consider the consequences of a weaker upper
bound on $V(1)$, such as $V(1) < 10 \omega_n$.  In this case, it
does not follow that $\phi_*([M]) = 0$ in $N$.  Instead, we get
an upper bound for the simplicial norm of $\phi_*([M])$ in
$H_n(N)$.

At this point, we recall the relevant facts about simplicial
norms.  For more information, see \cite{G5}.  Suppose that $C$
is a rational k-cycle in $M$.  We can write $C$ as a finite sum
$\sum a_i \Delta_i$, where $a_i$ is a rational number and
$\Delta_i$ is a map from the k-simplex to $M$.  We say that the
size of $C$ is equal to the sum $\sum |a_i|$.  The size of $C$ is
just the number of simplices counted with multiplicity.  Then we
define the simplicial norm of a homology class $h \in H(M,
\mathbb{Q})$ to be the infimal size of any rational cycle $C$ in
the class $h$.  We will write the simplicial norm of $h$ as $\| h
\|$.  For a closed oriented manifold $M$, the simplicial volume
of $M$ is defined to be the simplicial norm of the fundamental
class $[M]$.

We will use two facts about the simplicial norm.  The first
fact is that it decreases under any mapping.  In other words,
if $\phi: M \rightarrow N$ is a continuous map between spaces,
then $\| \phi_*(h) \| \le \| h \|$.  This property follows
immediately from the definition.  The second fact is that the
simplicial volume of a closed oriented hyperbolic manifold is
bounded below by the volume of the manifold.  We state this
result as a theorem.

\begin{reftheorem} (Thurston, see \cite{G5}) Suppose that $(M,
hyp)$ is a closed hyperbolic n-manifold.  Then the simplicial
volume of $M$ is at least $c(n) Vol(M,hyp)$.
\end{reftheorem}

(In fact, the simplicial volume of $M$ is equal to $c(n)
Vol(M,hyp)$ for an appropriate constant $c(n)$, but we don't need
this fact.)

Now we suppose that $(M^n,g)$ is a closed orientable Riemannian
manifold.  Let $N$ be the rectangular nerve constructed from a
good cover of $(M, g)$, and let $\phi: M \rightarrow N$ be the
map to the rectangular nerve.

\begin{lemma} For each number $V_0 > 0$ and each dimension $n$,
there is a constant $C(V_0, n)$ so that the following estimate
holds.  If $(M,g)$ has $V(1) < V_0$, then the simplicial norm $\|
\phi_*([M]) \| \le C(V_0, n) \textrm{ Volume}(M,g)$.
\end{lemma}

\proof The proof of this lemma is a modification of the proof of
Lemma 8.  According to Lemma 5, we have the following bound for
the volume of $\phi(M)$ in various regions of $N$.

$$| \phi(M) \cap Star(F)| < C_1 V(1) r_1(F)^{n+1} e^{- \beta
d(F)}.$$

\noindent In this formula, $C_1$ and $\beta$ are dimensional
constants, and $F$ can be any face of the rectangular nerve $N$. 
We let $\epsilon = \epsilon(\beta, n)$ be the same constant as in
the proof of Lemma 8.  We now divide the faces of $N$ into thick
and thin faces as follows.  If $C_1 V_0 r_1(F) < \epsilon$, then
we say that $F$ is thin, and otherwise we say that $F$ is thick.

We begin with some simple estimates about the thick and thin
simplices.  If $F$ is thin, then any higher-dimensional face
containing $F$ in its boundary is also thin.  Also, the dimension
of a thick face is bounded by $d(V_0, n)$, a constant depending only
on $n$ and $V_0$.  This estimate on the dimension follows from
Lemma 3, because if $I_{(0,1)}(F)$ contains $d(F)$ indices, then
the corresponding $d(F)$ balls contain a common intersection.  If
$B_1$ is the smallest of these balls, Lemma 3 guarantees that
$r_1 < C [\log d(F)]^{-1}$.  But by definition of a thick face,
$r_1(F) \ge c \epsilon / V_0$.

As before, we let $D$ be the dimension of $N$, and we define
$z_D$ to be the n-cycle $\phi(M)$.  We will construct a sequence
of homologous n-cycles $z_D \sim z_{D-1} \sim ... \sim z_n$, with
$z_k$ lying in the k-skeleton of $N$.  For any thin face $F$, the
cycle $z_k$ will obey the same estimate as in Lemma 8.

$$|z_k \cap Star(F)| < 2 \epsilon r_1(F)^n e^{- \beta d(F)}.$$

The construction is similar to the one in Lemma 8, but there is
an extra wrinkle having to do with the thick simplices.
First, we show that $z_D$ obeys the estimate that we want. 
Because of the definition of thin faces and the estimate in
equation 1, $|z_D \cap Star(F)| < \epsilon r_1(F)^n e^{- \beta
d(F)}$ for each thin face $F$.  

Now we suppose we have constructed $z_k$ for some $k>n$, and we
describe the construction of $z_{k-1}$.  Pick a k-dimensional
face $F$.  If $F$ is thin, then we define $z_k'
\cap F$ to be a minimal cycle with boundary $\partial (z_k \cap
F) \subset \partial F$.  As in the proof of Lemma 8, $z_k' \cap
F$ lies within the $s(F)$-neighborhood of $\partial F$, where
$s(F) = r_1(F) [2 \omega_n^{-1} \epsilon e^{-\beta k}]^{1/n}$. 
If $F$ is a thick face, then we define $z_k'$ by removing $z_k
\cap F$ and replacing it by a chain in $\partial F$ with the same
boundary as $z_k \cap F$.  According to a construction of
Federer-Fleming, we can choose a chain with volume bounded by
$G(V_0, n) |z_k \cap F|$.  (The constant in the Federer-Fleming
construction depends on the dimension $d(F)$.  We noted above
that for a thick face $d(F)$ is bounded by $d(V_0, n)$.  Also,
the Federer-Fleming construction gives a certain constant if we
do it in a cube.  We have to do it in a rectangular face.  If the
dimensions of the rectangular face are very uneven, the constant
can blow up.  In our case, though, $r_1$ is bounded below by a
constant depending only on $n$ and $V_0$, and all the dimensions
are bounded above by $(1/100)$.  Therefore, the stretching
factor $G(V_0, n)$ depends only on the dimension $n$ and
$V_0$.)

We then proceed as in the proof of Lemma 8.  We define
$\delta(k)$ to be $[2 \omega_n^{-1} \epsilon e^{- \beta
k}]^{1/n}$, and we define $z_{k-1}$ to be $R_{\delta(k)}(z_k')$. 
The cycle $z_{k-1}$ is homologous to $z_k$ and lies in the
(k-1)-skeleton of $N$.  By the same calculation as in Lemma 8, it
follows that $z_{k-1}$ obeys the volume estimate for thin faces
$F$: $|z_{k-1} \cap STAR(F)| < 2 \epsilon r_1(F)^n e^{- \beta
d(F)}$.

To get at the simplicial norm, we consider the cycle $z_n$.  It
lies in the n-skeleton of $N$.  The cycle $z_n$ is homologous to
a sum of n-faces of $F$, $\sum_i c_i F_i$, where $|c_i| < |z_n
\cap F_i|/|F_i|$.  Taking the barycentric triangulation of each
face, it follows that the simplicial norm of $z_n$ is bounded by $C
\sum_i |c_i|$.  If $F_i$ is a thin n-face, then the bound
$|z_n \cap STAR(F_i)| < 2 \epsilon r_1(F_i)^n$ guarantees that
$c_i = 0$.  For thick faces, $c_i$ may be non-zero.  If $F_i$ is
a thick n-face, then the volume $|F_i|$ is bounded below, and so
it follows that $\sum |c_i| < C(V_0, n) |z_n|$.

To bound $|z_n|$, we consider the increase of volume $|z_k| /
|z_{k+1}|$.  We form $z_k$ from $z_{k+1}$ by a surgery in each
(k+1)-face - yielding $z_{k+1}'$, followed by applying
$R_{\delta(k)}$.  For $k > d(V_0, n)$, all the surgeries occur in thin
faces, and therefore each surgery decreases the volume of the
cycle.  The application of $R_{\delta(k)}$ increases volume by at
most a factor $[1 - 2 \delta(k)]^{-n}$.  Therefore, $|z_{d(V_0,n)}| <
\prod_{k=d(V_0,n)}^\infty [1 - 2
\delta(k)]^{-n} |z_D|$.  By the same calculation as in Lemma 8,
the last expression is bounded by $2 |z_D|$.

If $k < d(V_0, n)$, then we have the following much weaker volume
bound: $|z_{k-1}| \le G(V_0, n) [1 - 2 \delta(k)]^{-n} |z_k|$. 
Therefore, the volume $|z_n| < 4 G(V_0,n)^{d(V_0,n)} |z_D|$. 
Finally, by Lemma 4, the volume $|z_D|$ is bounded by $C(n)
|(M,g)|$.  Assembling all the inequalities, we see that the
simplicial norm of $\phi_*([M])$ is bounded by $C(V_0, n)
|(M,g)|$. \endproof

Now we have enough ammunition to prove our second theorem.

\begin{theorem} For each dimension n, there is a number
$\delta(n) > 0$ so that the following estimate holds.  Suppose
that $(M^n, hyp)$ is a closed hyperbolic n-manifold and that $g$
is another metric on $M$, and suppose that $Vol(M,g) < \delta(n)
Vol(M,hyp)$.  Let $(\tilde M, \tilde g)$ denote the universal
cover of $M$ with the metric induced from $g$.  Then the
following inequality holds.

$$V_{(\tilde M, \tilde g)} (1) > V_{\mathbb{H}^n}(1).$$

\end{theorem}

\proof First we consider the special case that $M$ is oriented
and that each non-contractible curve in $(M,g)$ has length at least 1. 
In this case, $V_{(\tilde M, \tilde g)} (1) = V_{(M,g)}(1)$.  We
will assume that $V_{(M,g)}(1)$ is at most $V_{\mathbb{H}^n}(1) =
V_0$, and we need to prove that $Vol(M,g) \ge \delta(n)
Vol(M,hyp)$.

As usual, we choose a good cover of $(M,g)$.  We let $N$ be the
rectangular nerve of the cover and let $\phi: M \rightarrow N$ be
the map to the rectangular nerve constructed in section 3. 
According to Lemma 9, $C(n, V_0) Vol(M,g) \ge \| \phi_*([M]) \|$. 
We have assumed that the shortest non-contractible curve in
$(M,g)$ has length at least $1$.  According to Lemma 7, there is
a map $\psi: N \rightarrow M$ so that $\psi \circ \phi: M
\rightarrow M$ is homotopic to the identity.  In particular,
$\psi_*( \phi_*([M])) = [M]$.  Since the simplicial volume
decreases under maps, it follows that $\| \phi_*([M]) \|$ is
equal to $\| [M] \|$, the simplicial volume of $M$.  Finally,
since $M$ is closed and oriented, Thurston's theorem guarantees
that the simplicial volume of $M$ is at least $c(n) Vol(M,hyp)$. 
Putting together these inequalities, we see that $Vol(M,g) \ge
C(n, V_0)^{-1} c(n) Vol(M,hyp)$.  Since $V_0$, the volume of the
unit ball in $\mathbb{H}^n$, is itself a dimensional constant, we
see that $Vol(M,g) \ge \delta(n) Vol(M,hyp)$ as desired.

Next we consider the general case, with no restriction on the
lengths of non-contractible curves in $(M,g)$.  Again, we assume
that $V_{(\tilde M, \tilde g)}(1) \le V_0$, and we have to prove
that $Vol(M,g) \ge \delta(n) Vol(M, hyp)$.  Since $M$ admits a
hyperbolic metric, the fundamental group of $M$ is residually
finite.  (The group of isometries of hyperbolic n-space is a 
subgroup of $SL(N, \mathbb{C})$ for sufficiently large $N$, and any
finitely generated subgroup of $SL(N, \mathbb{C})$ is residually
finite according to \cite{M}.)  Therefore, we can choose a finite
cover $(\hat M, \hat g)$ so that $\hat M$ is oriented and so that
every non-contractible closed curve in $\hat M$ has length at
least 1.  Let $\hat hyp$ be the pullback of the hyperbolic metric
on $M$ to $\hat M$.  By assumption, $V_{(\tilde M, \tilde g)} (1)
\le V_0$.  Since the universal cover of $(\hat M, \hat g)$ is the
same as that of $(M,g)$, it follows that $V_{(\hat M, \hat g)}(1)
\le V_0$.  Now by the first case, we can conclude that $Vol(\hat
M, \hat g) \ge \delta(n) Vol(\hat M, \hat hyp)$.

Now if the covering map, $\pi: \hat M \rightarrow M$ has degree
$D$, then $Vol (M, g) = (1/D) Vol (\hat M, \hat g)$, and $Vol
(M,hyp) = (1/D) Vol (\hat M, \hat hyp)$.  Therefore, the last
inequality implies that $Vol (M,g) \ge \delta(n) Vol(M,hyp)$.
\endproof

Remark: There is a slightly more general result that holds with
the same proof.  It applies to products of hyperbolic manifolds. 
Suppose that $M^n$ is a product of closed manifolds, $M = M_1
\times ... \times M_d$, and that each manifold $M_i$ admits a
hyperbolic metric $hyp_i$.  Let $prod$ denote the product metric
$hyp_1 \times ... \times hyp_d$.  If $(M, g)$ has volume less
than $\delta(n) Vol(M, prod)$, then $V_{(\tilde M, \tilde g)}(1)
\ge V_{(\tilde M, \tilde prod)}(1)$.

\section{Open manifolds}

So far, we have proved Theorem 1 for closed manifolds.  Theorem 1
also holds for all complete Riemannian manifolds.  In this
section, we deal with the general case.  It requires only minor
technical modifications from the closed case.  We encourage the
reader not to take this section too seriously.

First, we review the definition of the filling radius of a
complete manifold.  The original definition appears on page 41
of \cite{G2}.  Let $(M^n, g)$ be a complete Riemannian
manifold.  The Kuratowski embedding maps a point $x \in M$ to the function
$dist_x$.  Since $M$ may not be compact, this function is unbounded.  
Nevertheless, it defines a measurable function, and the triangle inequality
implies that $|dist_x - dist_y|_{\infty} = dist(x,y)$.  The image of the Kuratowski
embedding lies in an affine copy of $L^\infty(M)$, namely all functions of the
form $dist_x + f$, where $f \in L^{\infty}(M)$.

Since $(M,g)$ is complete, any ball of finite radius is compact. 
The Kuratowski embedding is an isometry, and so the
inverse image of any compact set lies in a ball of finite
radius and is compact.  In other words, the Kuratowski
embedding is proper.  Therefore, the image of $M$ defines a cycle in
the sense of locally-finite homology theory.  The filling radius
of $(M,g)$ is the infimal $R$ so that this cycle bounds a locally
finite chain inside its $R$-neighborhood.  (If the cycle does not
bound within its $R$-neighborhood for every finite R, then the
filling radius is infinite.)

Most of the lemmas apply smoothly to complete manifolds with this definition, 
but a couple of them require some minor discussion.

Lemma 1 is local and applies immediately on a complete manifold.

Lemma 2 also holds on a complete manifold, but the proof requires a minor trick.
Let $K_1 \subset K_2 \subset ...$ be an exhaustion of $M$.  Consider the set of
balls $B(p, R/6)$ with $p \in K_i$ and $B(p,R)$ a good cover.  These balls cover $K_i$, and we can
find a finite subset of them that covers $K_i$.  Applying the Vitali covering
lemma, we get a cover of $K_i$ with all desired properties.  Repeating this procedure
for each i, we get a sequence of covers including more and more of $M$.  If we restrict attention
to the balls meeting a given compact set $K \subset M$ the set of possible covers is compact, which
we can check as follows.  Clearly, the set of possible centers $p \in K$ is compact.  We can
find a radius $r$ so that any ball of radius less than $r$ centered in $K$ has volume at least
$(1/2) \omega_n r^n$.  Therefore, every good ball has radius at least $r$ and at most $1/100$.
Also, every good ball has volume at least $(1/2) \omega_n r^n$.  Therefore, the number of balls
meeting $K$ is bounded.  We can therefore choose a subsequence so that the ball coverings restricted
to $K$ converge.  Now we again consider an exhaustion of $M$ by compact sets and diagonalize to
give a sequence of good coverings of $K_i$ that converge on all of $M$.  Their limit is our
good covering.

Lemma 3 is local and applies immediately on a complete manifold.

Lemma 4 is local as long as $U$ is bounded, and in this case it applies on a complete
manifold.

Lemma 5 follows immediately from Lemma 4.  It also holds on a complete manifold, except
possibly for the last estimate of the volume of $\phi(M)$.  This estimate
is not used in the proof of Theorem 1 anyway.

Lemma 6 follows for complete manifolds with the same proof.

Lemma 7 is not part of the proof of Theorem 1.

The most annoying technical problem occurs in Lemma 8.  The problem occurs
because the nerve $N$ may contain rectangles of every dimension.  Any given point will lie in
only finitely many balls, but as the point goes to infinity this
number may blow up.

In the original proof of Lemma 8, we had a cycle $z$ in $N$, and
we built a sequence of cycles $z = z_D \sim z_{D-1} \sim ... \sim
z_n$, where $D$ was the dimension of $N$ and $z_k$ lay in the
k-skeleton of $N$.  In general, the dimension of $N$ is not
finite, but is only locally finite, and we must proceed a little
differently.  Instead, we construct an infinite sequence of
cycles, $... \sim z_{k+1} \sim z_k \sim z_{k-1} \sim ... \sim
z_n$, with $z_k$ lying in the k-skeleton of $N$, so that the
sequence $z_k$ converges to $z$ as $k$ tends to infinity.

In a region of $N$ where the dimension is less than $k$, we
define $z_k$ to be the infinite composition $R_{\delta(k+1)}
\circ R_{\delta(k+2)} \circ ... $ applied to $z$.  (This infinite
composition is defined to be the limit of the maps
$R_{\delta(k+1)} \circ ... \circ R_{\delta(N)}$ as $N$ goes to
infinity.  The sequence of maps converges uniformly on compact
sets.)  In a region where the dimension of $N$ is at least $k$,
we define $z_k$ from $z_{k+1}$ as in the proof of Lemma 8.  Every
cycle $z_k$ is homologous to $z$ by a locally finite chain in
$N$, and the rest of the argument in Lemma 8 applies as before.

With these modifications, we get the proof of Theorem 1 in the
general case.

\begin{theorem1} (General case) For each dimension n, there is a
number $\delta(n) > 0$ so that the following estimate holds.  If
$(M^n,g)$ is a complete Riemannian n-manifold with filling radius
at least $R$, then $V(R) \ge \delta(n) R^n$.
\end{theorem1}

Finally, we prove a corollary about universal covers.

\begin{corollary} Let $(M^n, g)$ be a closed aspherical
Riemannian manifold, and let $V(R)$ measure the volumes of balls
in the universal cover $(\tilde M, \tilde g)$.  Then $V(R) \ge
\delta(n) R^n$ for all $R$.
\end{corollary}

\proof In \cite{G2} (page 43), Gromov proved that the
universal cover of $M$ has infinite filling radius.  Applying
Theorem 1, we get the corollary.
\endproof

\section{A question about Uryson width}

We say that the Ursyon k-width of a metric space $X$ is at most
$W$ if there is a continuous map $\pi$ from $X$ to a
k-dimensional polyhedron whose fibers have diameter at most $W$. 
The Uryson width is another way of measuring how ``thick'' a
manifold is, in a similar spirit to the filling radius.  Gromov
proved in \cite{G2} that Uryson (n-1)-width of a Riemannian
manifold $(M^n, g)$ controls its filling radius.  The opposite
inequality is not true, but it holds in many examples.  In other
words, an upper bound on the Uryson (n-1)-width is slightly
stronger than an upper bound on the filling radius.  It is an
open problem to understand how the volume of a Riemannian
manifold constrains its Uryson width.

\begin{ques} Is there a dimensional constant $C(n)$ so that every
closed Riemannian manifold $(M^n, g)$ has Uryson (n-1)-width at
most $C(n) Vol(M,g)^{1/n}$?
\end{ques}

This question is analogous to Gromov's estimate for the filling
radius.  There is another question, analogous to Theorem 1.

\begin{ques} Is there a dimensional constant $c(n)$ so that every
closed Riemannian manifold $(M^n,g)$ with Uryson (n-1)-width at least $W$
has $V(W) \ge c(n) W^n$?
\end{ques}

An affirmative answer to the second question is stronger than an
affirmative answer to the first question.

It looks plausibe that our proof of Theorem 1 can be modified to
bound the Uryson (n-1)-width of $M$ instead of its filling
radius, giving an affirmative answer to the second question.  To
bound the Uryson width, we would modify the proof of Lemma 8. For
each k-face $F$ of $N$, we make $z_k' \cap F$ the image of $z_k$
under a MAP that fixes the boundary $z_k \cap \partial F$ and
minimizes volume subject to the boundary restriction.  I believe
that $z_k'$ should be a minimal cycle plus a measure 0 region
that can be pushed as close to $\partial F$ as the minimal cycle
piece.  Then our argument gives a sequence of homotopies of the
map $\phi$, ending with a map $\pi$ into the (n-1)-skeleton of
$N$.  If $\pi(x)$ belongs to a face $F$, then $\phi(x)$ must have
belonged to $Star(F)$.  Therefore, if $I_+(F)$ contains an index
$i$, then $\pi^{-1}(F) \subset \phi^{-1}[Star(F)]$ lies in $B_i$. 
Since each ball in our cover has radius at most $(1/100)$, it
would follow that the Uryson (n-1)-width of $(M,g)$ is at most
$(2/100)$.

\end{document}